\documentclass[oneside]{amsart}
\usepackage[margin=35mm]{geometry}
\usepackage{subfigure}

\newif\ifhavetikz
\usepackage{tikz}
\usetikzlibrary{decorations.pathmorphing}
\usetikzlibrary{decorations.pathreplacing}
\usetikzlibrary{shapes.geometric}
\usetikzlibrary{arrows}
\pgfrealjobname{hompdo}
\havetikztrue

\usepackage{hyperref}
\hypersetup{pdfauthor={Dan Drake},
            pdftitle={Higher order matching polynomials and d-orthogonality}}

\newtheorem{theorem}{Theorem}[section]
\newtheorem{corollary}[theorem]{Corollary}
\newtheorem{proposition}[theorem]{Proposition}

\newcommand{\nonthmref}[2][Corollary]{\hyperref[#2]{#1~\ref*{#2}}}

\edef\Lslash{\L}
\renewcommand{\L}{\mathcal{L}}
\newcommand{\Lt}{\mathcal{L}\sbspt{}}

\newcommand{\luk}{\Lslash ukasiewicz}
\newcommand{\sbsp}[2]{_{#1}^{(#2)}}
\newcommand{\sbspt}[1]{\sbsp{#1}{t}}
\newcommand{\sbspnt}{\sbspt{n}}

\DeclareMathOperator{\cyc}{Cyc}

\newcommand{\R}{\mathbb{R}}
\newcommand{\dx}{\,\mathrm{d}x}

\ifhavetikz
  \tikzstyle{vertex}=[circle, draw=black, fill=black, inner sep=1]
  \tikzstyle{unvertex}=[circle, inner sep=1]
  \tikzstyle{edgestyle}=[]
  \newcommand{\drawedge}[3][]{\draw (#2) to [out=75, in=105,
    looseness=1.5, edgestyle, #1] (#3);}

\fi

\begin{document}

\title{Higher order matching polynomials and $d$-orthogonality}
\author{Dan Drake}
\address{Department of Mathematical Sciences\\
  Korea Advanced Institute of Science and Technology\\
  Daejeon, Korea}
\email{ddrake@member.ams.org}
\urladdr{http;//mathsci.kaist.ac.kr/~drake}
\keywords{matching polynomials, orthogonal polynomials, $d$-orthogonality}
\subjclass[2000]{Primary: 05E35; Seconday: 05C70, 33C45}
\date{16 November 2009}

\begin{abstract}
  We show combinatorially that the higher-order matching polynomials of
  several families of graphs are $d$-orthogonal polynomials. The
  matching polynomial of a graph is a generating function for coverings
  of a graph by disjoint edges; the higher-order matching polynomial
  corresponds to coverings by paths. Several families of classical
  orthogonal polynomials---the Chebyshev, Hermite, and Laguerre
  polynomials---can be interpreted as matching polynomials of paths,
  cycles, complete graphs, and complete bipartite graphs. The notion of
  $d$-orthogonality is a generalization of the usual idea of
  orthogonality for polynomials and we use sign-reversing involutions to
  show that the higher-order Chebyshev (first and second kinds),
  Hermite, and Laguerre polynomials are $d$-orthogonal. We also
  investigate the moments and find generating functions of those
  polynomials.
\end{abstract}

\maketitle
\section{Introduction and background}
\label{sec:intro}

A \emph{matching} of a graph is a subset of mutually disjoint edges in
the graph. Given a matching, we can assign a weight to each matching by
giving a weight of $-1$ to each edge in the matching and weight $x$ to
each vertex not adjacent to an edge in the matching, then multiplying
those weights together. We define the \emph{matching polynomial} of a
graph as the sum of weights of all matchings of the graph.

Matching polynomials have long been an object of interest in graph
theory, and it is well-known that the matching polynomials for some
classes of graphs---namely paths, cycles, complete graphs, and complete
bipartite graphs---are in fact classical orthogonal polynomials:
respectively, the Chebyshev polynomials of the second and first kinds,
Hermite polynomials, and Laguerre polynomials. See
\cite{sainte-catherine.viennot:combinatorial, godsil:hermite,
  godsil:theory, strehl:fibonacci} and also \cite[\S
4]{stanton:orthogonal}, which all treat the links between classical
orthogonal polynomials and matching polynomials.

The number of matchings of a graph was used by Hosoya to develop his
``topological index'' $Z$, which relates chemical properties of
hydrocarbons with their molecular structure. Later, Randi\'c, Morales,
and Araujo \cite{randic.morales.ea:higher-order} generalized the $Z$
index to the so-called higher-order Hosoya numbers by considering
coverings of graphs not by disjoint edges (which can be thought of as
paths of length one), but by paths of length two, three, and so on.
Araujo, Estrada, Morales, and Rada, starting from the higher-order
Hosoya numbers and working with Farrell's $F$-coverings
\cite{farrell:general, farrell:decomposition, farrell:path}, described
the higher-order matching polynomial of a graph, derived recurrence
relations, exact formulas, and also found expressions for those
polynomials as hypergeometric series
\cite{araujo.estrada.ea:higher-order}.

Apart from combinatorics and graph theory, Van Iseghem
\cite{iseghem:approximants} and Maroni \cite{maroni:lorthogonalite}
introduced a generalization of orthogonality for polynomials. A set of
polynomials $\{P_{n}\}_{n \ge 0}$ that is orthogonal in the usual sense
has an associated positive measure $\mu$ such that
\[
\int P_{n} P_{m}\,\mathrm{d}\mu = 0
\quad\text{if $n > m$, and}\quad
\int P_{n}^{2}\,\mathrm{d}\mu \ne 0.
\]
(We also demand that the degree of $P_{n}$ be $n$.) We call the integral
of $P_{n}^{2}$ the $L^{2}$ norm of $P_{n}$. For our purposes, instead of
providing a measure, it is equivalent to give a sequence of moments
$\{\mu_{n}\}_{n \ge 0}$ and define a linear functional $\L$ on the space
of polynomials by declaring $\L(x^{n}) = \mu_{n}$. A set of orthogonal
polynomials must satisfy a recurrence relation of the form
\begin{equation}
  \label{eq:general-op-recur}
  P_{n+1} = (x - b_{n}) P_{n} - \lambda_{n} P_{n-1}.
\end{equation}
Van Iseghem and Maroni defined the concept of $d$-orthogonality, or
orthogonality of dimension~$d$. Here $d$ is a positive integer, and we
say that a monic set of polynomials $\{P_{n}\}_{n \ge 0}$ is
$d$-orthogonal if there is a measure $\mu$ (or, for us, a sequence of
moments) such that
\begin{equation}
  \label{eq:d-op-integral}
  \int P_{n} P_{m}\,\mathrm{d}\mu = 0
  \quad\text{if $n > dm$, and}\quad
  \int P_{dn} P_{n}\,\mathrm{d}\mu = 1.
\end{equation}
Observe that usual orthogonal polynomials correspond to $d=1$. We will
commit a minor abuse of language and call the integral of $P_{dn}P_{n}$
the $L^{2}$ norm of $P_{n}$. Sets of $d$-orthogonal polynomials satisfy
a recurrence relation of order $d+2$ analogous to the one above.

In this paper, we establish the $t$-orthogonality for higher-order
matching polynomials corresponding to coverings of paths, cycles,
complete graphs, and complete bipartite graphs by paths with $t$ edges.
We will find formulas and combinatorial descriptions for the moments,
the recurrence relation, a sign-reversing involution that proves the
orthogonality and $L^{2}$ norms, and generating functions for the
moments and polynomials.

\subsection{Notation and terminology}
\label{sec:notation-terminology}

A path with $t$ edges and $t+1$ vertices will be called a ``$t$-path''.
Vertices of a graph not adjacent to an edge in the matching will often
be called ``fixed points''; the term comes from thinking of a matching
of a graph as giving an involution on the vertices. The set of integers
from $1$ to $n$, inclusive, will be written as ``$[n]$''. Finally, we
use $A \sqcup B$ for the disjoint union of two sets (usually graphs).

The vertices of our graphs are all labeled $1$ to $n$, and we often draw
matchings by arranging the vertices horizontally and drawing arcs for
edges or paths in the matching; we say that a matching is noncrossing
when such a diagram has no crossings. We will draw set partitions in a
similar manner, as Kasraoui and Zeng do in
\cite{kasraoui.zeng:distribution}.

\section{Warmup: Chebyshev polynomials of the second kind}
\label{sec:warmup-cheby2nd}

We begin with the Chebyshev polynomials of the second kind, which is
the simplest example. The basic combinatorics of Chebyshev polynomials
of the second kind were described by de~Sainte-Catherine and Viennot in
\cite{sainte-catherine.viennot:combinatorial, viennot:theorie,
viennot:combinatorial} and we review the theory to familiarize the
reader with our basic strategy and aims.

The Chebyshev polynomial of the second kind $U_{n}(x)$ is defined here
as the matching polynomial of a path with $n$ vertices. With this
normalization, they are also called Fibonacci polynomials, since
matchings of a path with $n$ vertices corresponds in a natural way to
``pavings'' of length $n$ composed of dominos and monominos, and such
pavings are counted by Fibonacci numbers. The last vertex of such a
matching must be a fixed point or in an edge in the matching, so the
recurrence relation corresponding to \eqref{eq:general-op-recur} is
clear:
\begin{equation}
  \label{eq:cheby2nd-recurrence}
  U_{n+1}(x) = x U_{n}(x) - U_{n-1}(x),
\end{equation}
so these polynomials have $b_{n} = 0$ and $\lambda_{n} = 1$ for all $n$.
Viennot established that the $n$th moment of a set of orthogonal
polynomials with recurrence coefficients $b_{n}$ and $\lambda_{n}$ as in
\eqref{eq:general-op-recur} equals the total weight of all weighted
Motzkin paths of length $n$; a Motzkin path is a lattice path that never
goes below the $x$-axis and takes upsteps, horizontal steps, and
downsteps (that is, steps of the form $(1,1)$, $(1,0)$, and $(1,-1)$),
with upsteps of weight $1$, horizontal steps at height $n$ of weight
$b_{n}$, and downsteps leaving from height $n$ of weight $\lambda_{n}$.
Knowing that, we see that the $n$th moment of the Chebyshev polynomials
of the second kind is the number of Dyck paths---lattice paths with only
up- and down-steps---of length $n$.

The number of Dyck paths of length $2m$ is the Catalan number
$\binom{2m}{m}/(m+1)$, which here we interpret as the number of
noncrossing complete matchings of $K_{2m}$, which are the same as
noncrossing set partitions of $[2m]$ in which all blocks have size two.

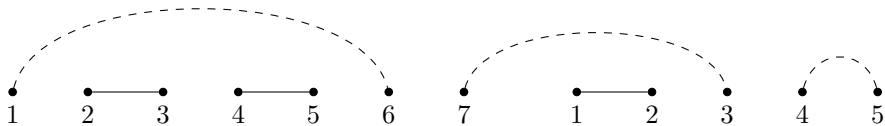
\begin{figure}[ht]
\centering
\beginpgfgraphicnamed{tikz/cheby2nd-integral-out}%
\begin{tikzpicture}
  \foreach \a in {1,...,7}{
    \node (L\a) at (\a,0) [vertex, label=below:$\a$] {} ;
  }

  \foreach \a in {1,...,5}{
     \node (R\a) at (\a + 7.5, 0) [vertex, label=below:$\a$] {} ;
  }

  \drawedge[dashed, looseness=.75]{L1}{L6}
  \draw (L2) -- (L3) ;
  \draw (L4) -- (L5) ;
  \drawedge[dashed, looseness=.75]{L7}{R3}
  \draw (R1) -- (R2) ;
  \drawedge[dashed]{R4}{R5}
\end{tikzpicture}

\endpgfgraphicnamed
\caption{A configuration that contributes weight $-1$ to the
  ``integral'' $\L(U_{7}(x) U_{5}(x))$. The orthogonality involution
  would change the edge connecting vertices $2$ and $3$ on the left from
  solid to dashed. There must always be at least one homogeneous
  adjacent edge, so all configurations are canceled and the integral is
  zero.}
\label{fig:cheby2nd-integral}
\end{figure}

The orthogonality of these Chebyshev polynomials can be proved with a
sign-reversing involution: $U_{n}(x) U_{m}(x)$ is the generating
function for pairs of matchings of an $n$-vertex path and and $m$-vertex
path. Integrating that product can be interpreted as the generating
function for complete noncrossing matchings on $[n] \sqcup [m]$ with:
\begin{itemize}
\item black edges of weight $-1$ that connect adjacent vertices and are
  homogeneous---that is, they stay within $[n]$ or $[m]$, and
\item dashed edges of weight $1$ between any two vertices.
\end{itemize}
The weight of such a configuration is the product of the weights of the
edges. See \autoref{fig:cheby2nd-integral} for an example. Given such a
configuration, we can produce another configuration by finding the
leftmost edge that connects adjacent vertices and changing it to black
if it is dashed, or vice versa. This process is a sign-reversing
involution that cancels all configurations with a homogeneous edge, so
$\L(U_{n}(x)U_{m}(x))$ equals the number of complete noncrossing
inhomogeneous matchings of $[n] \sqcup [m]$. If $n \ne m$, there are
obviously zero such matchings, and if $n = m$, there's exactly one: a
``rainbow'' configuration in which vertex $n-k$ on the left is connected
to vertex $k$ on the right ($1 \le k \le n$).

Let's finish this section by mentioning the generating functions of the
polynomials and the moments. For $U_{n}(x)$, any such polynomial is a
sequence of fixed points, which have size $1$ and weight $x$, and edges,
which have size $2$ and weight $-1$. When your objects are composed of
sequences of smaller objects, the generating function is typically just
a geometric series:
\begin{equation}
\label{eq:cheby2-poly-gf}
\sum_{n \ge 0} U_{n}(x) z^{n} = \frac{1}{1 - (xz - z^{2})}.
\end{equation}
The moments are Catalan numbers, whose generating function is well
known; for example, see Aigner \cite[\S 3.1 and \S 7.3]{aigner:course}:
\begin{equation}
\label{eq:cheby2-mu-gf}
f(z) = \sum_{n \ge 0} \mu_{n} z^{n} = \frac{1 - \sqrt{1 - 4 z^{2}}}{2 z^{2}},
\end{equation}
but for our purposes, we will focus more on the functional equation
satisfied by $f(z)$ and the corresponding continued fraction. The
functional equation is
\begin{equation}
\label{eq:cheby2-mu-functional-eq}
f(z) = 1 + z^{2} f(z)^{2}
\end{equation}
and is easy to explain using Dyck paths, which are counted by the
Catalan numbers. Think of $f(z)$ as standing for ``any possible Dyck
path''; such a path is either empty (with weight $1$), or is of the form
``upstep-(some Dyck path)-downstep-(some Dyck path)'', which has weight
$z^{2} f(z)^{2}$; see \autoref{fig:dyck-path-decomp}.

\newcommand{\upstep}[1]{\draw (#1) node[vertex] {} --
  node [midway, above] {$z$} ++(1,1) node[vertex] {};}

\newcommand{\wiggle}[3]{\draw [decorate, decoration=snake] (#1) to
  [out=45, in=135] node [midway, above, outer sep=1mm] {$#3$} (#2);}

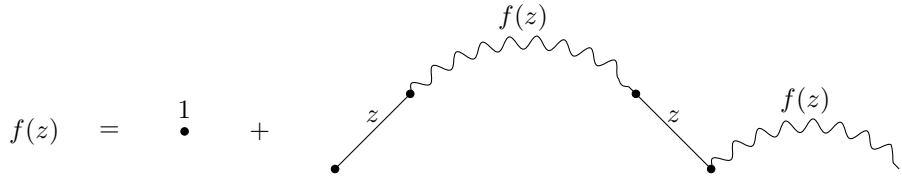
\begin{figure}[ht]
  \centering \beginpgfgraphicnamed{tikz/dyck-path-decomp-out}%
\begin{tikzpicture}

\node at (-4, 1/2) {$f(z)$};

\node at (-3, 1/2) {$=$};

\node at (-2, 1/2) [vertex, label=above:$1$] {} ;

\node at (-1, 1/2) {$+$} ;

\upstep{0,0}
\wiggle{1, 1}{4, 1}{f(z)}
\draw (4, 1) node[vertex] {} -- node [midway, above] {$z$} ++(1, -1)
  node[vertex] {};
\wiggle{5, 0}{7.5, 0}{f(z)}

\end{tikzpicture}
\endpgfgraphicnamed
  \caption{A pictorial explanation of the functional equation of
    \eqref{eq:cheby2-mu-functional-eq}.}
  \label{fig:dyck-path-decomp}
\end{figure}

By rearranging \eqref{eq:cheby2-mu-functional-eq}, one is easily led to
a continued fraction expression for $f(z)$:
\begin{equation}
  \label{eq:cheby2-mu-contfrac-gf}
  f(z) = \cfrac{1}
  {1 - \cfrac{z^{2}}
    {1 - \cfrac{z^{2}}
      {1 - \cdots\vphantom{\cfrac{z^{2}}{1}}}}}
\end{equation}
Now we generalize the above work to higher-order matching polynomials.

\section{Higher-order Chebyshev polynomials of the second kind}
\label{sec:high-order-cheby2nd}

Following Araujo et al. in \cite{araujo.estrada.ea:higher-order}, let's
now cover the path with $n$ vertices by paths with $t$ edges, and give
$t$-paths weight $-1$ and fixed points weight $x$. We will denote the
generating function for such coverings by $U_{n}^{(t)}(x)$ and call them
\emph{Chebyshev polynomials of the second kind and order~$t$}.

These polynomials satisfy a recurrence relation similar to
\eqref{eq:cheby2nd-recurrence}:
\begin{equation}
  \label{eq:cheby2-t-recurrence}
  U\sbspt{n+1}(x) = x U\sbspnt(x) - U\sbspt{n-t}(x);
\end{equation}
the proof is effectively the same: vertex $n+1$ is either fixed or the
final vertex in a $t$-path, and the rest of the vertices can be
covered by a smaller configuration.

Let $\mu\sbspnt$ be the number of noncrossing set partitions of $[n]$ in
which all blocks have size $t+1$, and let $\Lt$ be the linear functional
on the space of polynomials defined by $\Lt(x^{n}) = \mu\sbspnt$. Then
we have our first theorem, which generalizes a result of
de~Sainte-Catherine and Viennot \cite[Theorem~7]
{sainte-catherine.viennot:combinatorial}.

\begin{theorem}
  \label{thm:cheby2-t-integral}
  Let $n_{1}, n_{2},\dots,n_{k}$ be nonnegative integers. The integral
  \[
  \Lt\left(\prod_{i=1}^{k} U\sbspt{n_{i}} (x)\right)
  \]
  equals the number of inhomogeneous noncrossing coverings of $[n_{1}]
  \sqcup \dots\sqcup [n_{k}]$ by $t$-paths, or equivalently, the number
  of noncrossing set partitions of $[n_{1}] \sqcup\dots\sqcup [n_{k}]$
  in which all blocks have size $t+1$ and no block is a subset of any
  $[n_{i}]$.
\end{theorem}

\begin{proof}
  The proof uses a sign-reversing involution analogous to the above
  involution we used for the usual Chebyshev polynomials of the second
  kind. The product of the polynomials is the generating function for
  $k$-tuples of coverings of $n_{i}$-vertex paths by fixed points and
  $t$-paths. Integrating the product yields the generating function for
  complete noncrossing coverings of $[n_{1}] \sqcup\dots\sqcup [n_{k}]$
  by $t$-paths with
  \begin{itemize}
  \item homogeneous black $t$-paths of weight~$-1$ that connect a
    sequence of adjacent vertices in the underlying path, and
  \item dashed $t$-paths of weight~$1$ that can go anywhere.
  \end{itemize}
  Call a path (black or dashed) that connects $t+1$ adjacent vertices
  in the underlying path ``flat''. The sign-reversing involution is
  simple: find the leftmost homogeneous flat $t$-path and change its
  ``color'' from black to dashed, or vice versa. Any configuration
  that has at least one homogeneous $t$-path must have a flat
  $t$-path, so this involution will cancel any configuration with a
  homogeneous edge. Uncanceled configurations have only edges of
  weight~$1$, so $\Lt\big(U\sbspt{n} (x) U\sbspt{m} (x)\big)$ equals
  the number of configurations with only inhomogeneous edges.
\end{proof}

The above theorem immediately implies that the polynomials $U\sbspnt$
are, in fact, $t$-orthogonal with respect to those moments.

\begin{corollary}
  \label{thm:cheby2-t-orthog}
  The polynomials $U\sbspnt$ are $t$-orthogonal with respect to the
  above moments: if $m > nt$, then
  \begin{equation}
    \label{eq:cheby2-t-orthog}
    \L\big(U\sbspt{m} (x) U\sbspt{n} (x)\big) = 0
    \quad \text{and} \quad
    \L\big(U\sbspt{nt} (x) U\sbspt{n} (x)\big) = 1.
  \end{equation}
\end{corollary}

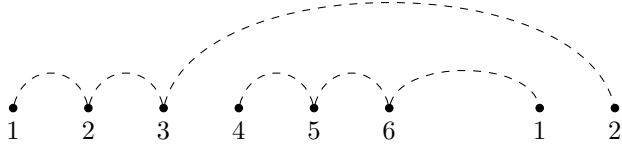
\begin{figure}[h]
  \centering
  \beginpgfgraphicnamed{tikz/cheby2nd-t-l2norm-out}%
\begin{tikzpicture}
  \foreach \a in {1,...,6}{
    \node (L\a) at (\a,0) [vertex, label=below:$\a$] {} ;
  }

  \foreach \a in {1,...,2}{
     \node (R\a) at (\a + 7, 0) [vertex, label=below:$\a$] {} ;
  }

  \drawedge[dashed]{L1}{L2}
  \drawedge[dashed]{L2}{L3}
  \drawedge[dashed, looseness=.8]{L3}{R2}
  \drawedge[dashed]{L4}{L5}
  \drawedge[dashed]{L5}{L6}
  \drawedge[dashed, looseness=.8]{L6}{R1}
\end{tikzpicture}
\endpgfgraphicnamed
  \caption{The sole uncanceled configuration in $\L \big(
      U\sbsp{6}{3}(x) U\sbsp{2}{3}(x) \big)$.}
  \label{fig:cheby2nd-t-l2norm}
\end{figure}

\begin{proof}
  If $m > tn$, then any configuration of black and dashed edges must
  have at least one homogeneous adjacent $t$-path in $[m]$, so the
  integral is zero. The integral of $U\sbspt{n} U\sbspt{nt}$ is~$1$
  because there is exactly one inhomogeneous configuration, an example
  of which is pictured in \autoref{fig:cheby2nd-t-l2norm}.
\end{proof}

The usual Chebyshev polynomials of the second kind are the generating
functions for matchings of a path; the number of such matchings is a
Fibonacci number, and can be obtained with the appropriate substitution:
$F_{n} = U_{n}(i)/i^{n}$. The higher-order polynomials lead to the
higher-order Fibonacci numbers of Randi\'c et al.
\cite{randic.morales.ea:higher-order}, and can be obtained by
\begin{equation}
  \label{eq:cheby2-t-fibonacci}
  {}^{t} F_{n} = \left(\frac{1}{w}\right)^{n} U\sbspnt (w),
\end{equation}
where $w$ is any $(t+1)$th root of $-1$ and ${}^{t} F_{n}$ is the
notation of Randi\'c et al. for the higher-order Fibonacci numbers. The
explanation for this is simple: we want each $t$-edge path to have
weight $+1$, but $U\sbspnt(x)$ gives them weight $-1$. By multiplying by
$1/w^{n}$, we effectively give each vertex weight $1/w$, and hence the
total weight of each $t$-path is $+1$. Then we plug in $w$ to give each
fixed point weight $+1$ as well.

Before we start investigating the moments of these polynomials, note
that the generating function for the polynomials is a straightforward
generalization of \eqref{eq:cheby2-poly-gf}.

\begin{proposition}
  \label{thm:cheby2-t-poly-gf}
  The ordinary generating function of the Chebyshev polynomials of the
  second kind and order $t$ is
  \begin{equation}
    \label{eq:cheby2-t-poly-gf}
    UP(t,x,z) = \sum_{n \geq 0} U\sbspnt (x) z^{n} =
    \frac{1}{1-(xz - z^{t+1})}.
  \end{equation}
\end{proposition}

The ``$UP$'' is intended to be mnemonic: the $U$ is for $U\sbspnt(x)$,
and the $P$ is for ``polynomials''. We'll meet $UM$, the generating
function for the moments, shortly.

\begin{proof}
  The rational function in \eqref{eq:cheby2-t-poly-gf} equals
  \[
  \sum_{k \geq 0} \left( xz - z^{t+1} \right)^{k},
  \]
  which can be interpreted as the generating function for finite
  sequences of objects with either weight $xz$ or weight $-z^{t+1}$
  grouped by number of objects; the sum in \eqref{eq:cheby2-t-poly-gf}
  is the generating function for the same thing, just grouped by
  coefficient of $z$.
\end{proof}

\subsection{Moments of \texorpdfstring{$U\sbspnt (x)$}{higher order
    U's}}
\label{sec:cheby2-moments}

The moments for the usual Chebyshev polynomials of the second kind
(order $t=1$) are ``aerated'' Catalan numbers: $\mu_{2n+1} = 0$ and
$\mu_{2n} = \binom{2n}{n}/(n+1)$ for all nonnegative $n$. In what
follows, we will work with the following interpretations of Catalan
numbers: as Dyck paths with $2n$ steps, as noncrossing matchings of
$[2n]$, as binary trees, and as triangulations of an $(n+2)$-gon.

To understand the moments for Chebyshev polynomials of the second kind
and order $t$, we need to understand the \emph{Fuss-Catalan numbers},
also called \emph{$k$-Catalan} or \emph{generalized Catalan numbers}.
The Fuss-Catalan numbers $C\sbsp{n}{k}$ are defined by
\begin{equation}
  \label{eq:fuss-catalan-num}
  C\sbsp{n}{k} = \frac{1}{kn + 1} \binom{(k+1)n}{n},
\end{equation} and have been extensively studied; see Hilton and
Pedersen \cite{hilton.pedersen:catalan} for an introduction to these
numbers, which were likely first described by Fuss \cite{fuss:solutio}
(see the table in that paper on page $249$), nearly $50$ years before
Catalan \cite{catalan:note}. The above interpretations of the usual
Catalan numbers generalize to the following interpretations for the
Fuss-Catalan numbers:
\begin{itemize}
\item Dissections of an $(nk+2)$-gon into $(k+2)$-gons.
\item Rooted plane trees in which all non-leaf vertices have $k+1$
  children; i.e., $(k+1)$-ary trees.
\item Lattice paths of length $n$ composed of $(1,1)$ and $(1,-k)$ steps
  that start and end on the $x$-axis, and never go below the $x$-axis.
\item Noncrossing set partitions with all blocks of size $k+1$.
\end{itemize}
The first item was studied by Fuss; the second and third are connected
by the \luk{} language; see \cite[chap. 11]{lothaire:combinatorics}. The
last item is not as ubiquitous as the others and we will use it to prove
the following.

\begin{theorem}
  \label{thm:cheby2-moments}
  The moments $\mu\sbspnt$ for $U\sbspnt (x)$ are ``aerated''
  Fuss-Catalan numbers: $\mu\sbspt{(t + 1)n} = C\sbspt{n}$ for all $n
  \ge 0$, and are zero otherwise.
\end{theorem}

\begin{figure}[ht]
  \centering
  \beginpgfgraphicnamed{tikz/ternary-tree-to-nc-setp-out}%
\begin{tikzpicture}[every node/.style={vertex}, level/.style={sibling
      distance=25mm/#1}, level distance=10mm]

  \node (root) {} {
     child {
        child
        child {
           child
           child
           child }
        child}
     child
     child {
        child
        child
        child }
     };

  \node at (root-1) [label=left:$1$] {};
  \node at (root-1-1) [label=below:$2$] {};
  \node at (root-1-2) [label={[inner sep=1mm]left:$3$}] {};
  \node at (root-1-2-1) [label=below:$4$] {};
  \node at (root-1-2-2) [label=below:$5$] {};
  \node at (root-1-2-3) [label=below:$6$] {};
  \node at (root-1-3) [label=below:$7$] {};
  \node at (root-2) [label=below:$8$] {};
  \node at (root-3) [label={[inner sep=1.5mm]left:$9$}] {};
  \node at (root-3-1) [label=below:$10$] {};
  \node at (root-3-2) [label=below:$11$] {};
  \node at (root-3-3) [label=below:$12$] {};
\end{tikzpicture}

\endpgfgraphicnamed
  \caption{An example illustrating the bijection from ternary trees to
    noncrossing set partitions in which all blocks have size three; the
    above tree corresponds to the set partition $\{\{1, 8, 9\},
    \{2,3,7\}, \{4, 5,6\}, \{10,11,12\}\}$.}
  \label{fig:ternary-tree-to-nc-setp}
\end{figure}
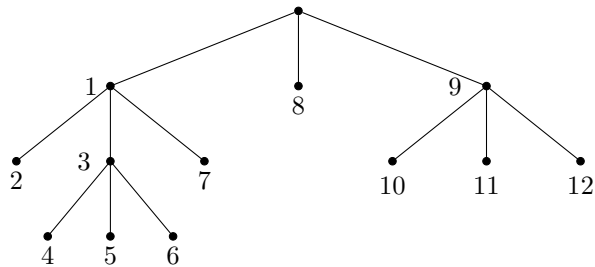

\begin{proof}
  The moments for a $d$-orthogonal set of polynomials are unique (this
  follows from an inductive argument, or by the ``spanning argument'' of
  \cite{roman.rota:umbral}), so from \autoref{thm:cheby2-t-orthog}, we
  know that noncrossing set partitions with all blocks of size $t+1$ are
  the correct moments for $U\sbspnt$. Therefore we need nothing more
  than a bijection to one of the above families of objects counted by
  the aerated Fuss-Catalan numbers. One easy bijection is to $(t+1)$-ary
  trees. Given a $(t+1)$-ary tree on $(t+1)n + 1$ vertices, number the
  vertices of the tree according to a depth-first, left-to-right search.
  The labels on each set of $t+1$ siblings in the tree describe the
  corresponding block in the set partition---see
  \autoref{fig:ternary-tree-to-nc-setp}. Injectivity is obvious, and
  surjectivity follows from the nesting structure of such a set
  partition. Consider the usual diagram of the set partition, and find
  all blocks of the set partition that contain $t+1$ consecutive
  numbers; those blocks will be sets of siblings in the tree that have
  no subtree. The root of the subtree corresponding to a block $P$ gets
  connected to vertex $i$ corresponding to block $Q$ when the smallest
  number in $P$ is one larger than the $i$th smallest number in $Q$.
  (For example, in \autoref{fig:ternary-tree-to-nc-setp}, the block
  $\{4,5,6\}$ is immediately nested by block $\{2,3,7\}$, and we connect
  the root of the $\{4,5,6\}$ subtree to $3$ because $4$ is one larger
  than $3$.) This produces a $(t+1)$-ary tree from a noncrossing set
  partition with all blocks of size $t+1$, so surjectivity---hence
  bijectivity---holds.
\end{proof}

There is a closely related bijection, analogous to Flajolet's path
diagrammes \cite{flajolet:combinatorial} and Viennot's Laguerre
histories \cite{viennot:theorie,viennot:combinatorial}.

A \luk{} path is a generalization of a Motzkin path in which, in
addition to upsteps and horizontal steps, downsteps of the form $(1,
-k)$ are allowed. The paths begin and end on the $x$-axis and never go
below the $x$-axis. For a general set of $d$-orthogonal polynomials, the
$n$th moment is the generating function for weighted \luk{} paths of
length $n$ that have downsteps $(1, -1), (1,-2), \dots, (1,-d)$. The
weight of each horizontal step and downstep is given by the recurrence
coefficients for the polynomials. See \cite[\S 4.2]
{roblet:interpretation} for more details.

The recurrence relation \eqref{eq:cheby2-t-recurrence} tells us that the
moments for Chebyshev polynomials of the second kind and order $t$ are
also given by \luk{} paths with upsteps $(1,1)$ and downsteps $(1, -t)$,
with all steps of weight $1$. The bijection between those paths and the
set partitions counted by $\mu\sbspnt$ is the obvious generalization of
the classical Motzkin path bijection to set partitions: given such a
\luk{} path with $n$ steps, one produces a set partition of $[n]$ with
the following procedure: begin with an empty set partition and read
through the path. if the $k$th step is an upstep, add $k$ to a set of
``candidates''. If step $k$ is a downstep, add a block to the set
partition consisting of $k$ and the $t$ largest elements of the
candidate set. This produces a set partition in which every block has
size $t+1$. Choosing the $t$ largest candidates guarantees that the set
partition is noncrossing, and the procedure is a bijection because every
downstep corresponds to a unique set of $t$ upsteps. This is in some
sense a generalization of the map between set partitions and Charlier
diagrams found in \cite[\S 3.1]{kasraoui.zeng:distribution} and
\cite{flajolet:combinatorial, viennot:theorie}; in the language of that
bijection, we place an ``opener'' vertex whenever one sees an upstep,
and placing a ``closer'' vertex and connecting the rightmost $t$ open
vertices.

For example, if $U$ stands for an upstep and $D$ for a $(1, -2)$ step,
the tree and set partition in \autoref{fig:ternary-tree-to-nc-setp}
correspond to the lattice path $UUUUUDDUDUUD$.

The tree and lattice path representations for the moments of $U\sbspnt$
make the generalizations of \eqref{eq:cheby2-mu-functional-eq} and
\eqref{eq:cheby2-mu-contfrac-gf} obvious.

%
%
%
\newsavebox{\fracbox}
\newlength{\fracraise}
\newlength{\myht}
\newlength{\mydp}

\newcommand{\mylr}[1]{%
\sbox{\fracbox}{\ensuremath{#1}}%
\settoheight{\myht}{\usebox{\fracbox}}%
\settodepth{\mydp}{\usebox{\fracbox}}%
\setlength{\fracraise}{\mydp}%
\addtolength{\fracraise}{-\myht}%
\setlength{\fracraise}{.5\fracraise}%
\addtolength{\fracraise}{.25em}%
\left(\raisebox{\fracraise}{\usebox{\fracbox}}\right)}

\begin{theorem}
  \label{thm:cheby2-t-moment-gf}
  Let $UM(t, z)$ be the ordinary generating function of the moments of
  the Chebyshev polynomials of the second kind and order $t$:
  \[
  UM(t, z) = \sum_{n \ge 0} \mu\sbspt{n} z^{n}.
  \]
  Then $UM(z)$ satisfies
  \begin{equation}
    \label{eq:cheby2-t-moment-gf-funceqn}
    UM(t, z) = 1 + z^{t+1} UM(t, z)^{t+1}
  \end{equation}
  and has the continued fraction expansion

  \begin{equation}
    \label{eq:cheby2-t-moment-gf-contfrac}
    UM(t, z) = \cfrac{1}{
      1 - \cfrac{z^{t+1}}{
        \mylr{1 - \cfrac{z^{t+1}}{
            \mylr{1 - \cfrac{z^{t+1}}{(1 - \cdots)^{t}}}^{t}}}^{t}}}.
  \end{equation}
\end{theorem}

Viennot called the above continued fraction an L-fraction when he
derived a generalization of the above theorem in \cite[chapter V, \S
6]{viennot:theorie}. Note that the continued fraction expansion only
requires the recursion coefficients of the polynomials, which Araujo et
al. found \cite[eq. (2.7)]{araujo.estrada.ea:higher-order}, and are
clear from the combinatorial description.

For these polynomials, and every class of higher-ordering matching
polynomials we consider here, we note that Araujo et al.
\cite{araujo.estrada.ea:higher-order} give exact formulas and
expressions for these polynomials as generalized hypergeometric
functions. They give an exact formula for $U\sbspnt(x)$ in equation
(2.11) and express that polynomial as a $_{t+1} F_{t}$ in (3.6).

\section{Chebyshev polynomials of the first kind}
\label{sec:cheby1-t}

Now let's move on to the Chebyshev polynomials of the first kind and
order $t$, which are the higher-order matching polynomials for a cycle
with, as usual, weight $x$ for fixed points and weight~$-1$ for a path
with $t$ edges. We denote them $T\sbspnt(x)$, and use $\cyc(n)$ for the
underlying labeled $n$-cycle. Here a $1$-cycle is a single
vertex\footnote{One can also have a $1$-cycle be a vertex with a loop,
with the convention that such a graph has no one-edge matching since the
edge would be incident with itself---but that's effectively the same as
just saying it's a single vertex with no edge at all.} and a $2$-cycle,
of course, is two vertices with two edges between them. The
interpretation of Chebyshev polynomials as the matching polynomial of a
cycle has been studied by several authors; see Benjamin and Walton
\cite{benjamin.walton:counting}, Bergeron \cite{bergeron:combinatoire},
and Hosoya \cite{hosoya:topological, hosoya:graphical}, whose $Q(Y)$
polynomials for cycloparaffins in the first reference are essentially
Chebyshev polynomials of the first kind. The Chebyshev polynomials of
the first kind are often defined by $P_{n}(\cos\theta) = \cos(n
\theta)$, and those polynomials are related to ours by $T\sbsp{n}{1}(x)
= 2 P_{n}(x/2)$.

The combinatorial model for these polynomials is well known, but it is
much harder to find models for the moments and an involution proof like
that of \autoref{thm:cheby2-t-integral} to prove orthogonality. The
Chebyshevs of the first kind satisfy the same recurrence relation as
those of the second kind, but they have different initial conditions.
Before we describe the moments and orthogonality involution, let's show
directly that $T\sbspnt (x)$ satisfies the same recurrence as
$U\sbspnt(x)$.

\subsection{A weight-preserving bijection for the recurrence relation}
\label{sec:cheby1-t-recurrence}

In this section, we find a weight-preserving bijection that shows
\begin{equation}
  \label{eq:cheby1-recurrence}
  T\sbspt{n+1}(x) = x T\sbspt{n}(x) - T\sbspt{n-t}(x),
\end{equation}
with $T\sbspt{i}(x) = x^{i}$ for $0 \le i \le t$ and $T\sbspt{t+1}(x) =
x^{t+1} - (t+1)$. Araujo et al. \cite[\S
2.2]{araujo.estrada.ea:higher-order} and Farrell \cite[\S
6]{farrell:path} both find this recurrence, but here we present a direct
bijection. The initial conditions are clear; either there are not enough
edges to have a $t$-path (so every vertex is fixed), or there are $t+1$
edges, and we can choose any one of them to be the single edge not in a
$t$-path.

Assume $n > t$. There are four cases for the bijection between coverings
of $T\sbspt{n+1}(x)$ and the union of coverings of $T\sbspnt(x)$ with
extra weight $x$ and coverings of $T\sbspt{n-t}(x)$ with extra weight
$-1$.

Given a configuration for $T\sbspnt(x)$, the edge from $1$ to $2$ is
either in a $t$-path or not. Case (a) is when the edge is not in a
$t$-path and is illustrated in \autoref{fig:cheby1-t-recurrence-a}. In
all four figures, we have a section of $\cyc(n+1)$ on top, and a section
of $\cyc(n)$ or $\cyc(n-t)$ on the bottom; the rest of the vertices are
omitted for clarity. Black edges are in a $t$-path, and dashed edges are
not. The gray triangles indicate expanding (or contracting) edges in the
map from one configuration to another. In case (a), we simply insert a
new edge ``behind'' vertex $1$ that is not in a $t$-path.

In case (b), the edge from $1$ to $2$ in $\cyc(n)$ is in a $t$-path,
which extends back to vertex $k$. We expand vertex $k-1$ into an edge
not in a $t$-path and relabel the relevant vertices in the $t$-path.
\autoref{fig:cheby1-t-recurrence-b} shows the process. Observe that this
always yields a configuration in which $(1,2)$ is in a $t$-path, and
that path is preceded by at least two edges not in a path.

Both operations multiply the weight by $x$ and yield a configuration for
$T\sbspt{n+1}(x)$.

\ifhavetikz
  \tikzstyle{contract}=[black!20]
  \tikzstyle{newvert}=[vertex, label={[font=\footnotesize, rectangle,
    inner sep=2pt]90:$#1$}]
  \tikzstyle{cheby1scale}=[scale=1.5]
\fi

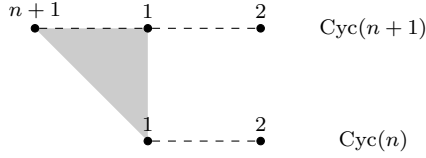
\begin{figure}
  \centering
  \beginpgfgraphicnamed{tikz/cheby1-recurrence-a-out}%
\begin{tikzpicture}[cheby1scale]
  \fill[contract] (0,0) -- (1,0) -- (1,-1);
  \node[font=\footnotesize] at (3, 0) {$\cyc(n+1)$};
  \node[font=\footnotesize] at (3, -1) {$\cyc(n)$};

  \foreach \v/\x in {n+1/0, 1/1, 2/2}
    \node[newvert=\v] (t\v) at (\x, 0) {};

  \foreach \v/\x in {1/1, 2/2}
    \node[newvert=\v] (b\v) at (\x, -1) {};

  \draw[dashed] (tn+1) -- (t1);
  \draw[dashed] (t1) -- (t2);

  \draw[dashed] (b1) -- (b2);
\end{tikzpicture}

\endpgfgraphicnamed
  \caption{Case (a) for the $T\sbspt{n+1}(x)$ recurrence. When the edge
    $(1,2)$ in $\cyc(n)$ is not covered by a $t$-path, we simply insert
    a new edge behind $1$, label the new vertex $n+1$, and attach the
    previous $(n,1)$ edge from $n$ to $n+1$. The weight of the new
    configuration is $x$ times the weight of the old.}
  \label{fig:cheby1-t-recurrence-a}
\end{figure}

\begin{figure}
  \centering
  \beginpgfgraphicnamed{tikz/cheby1-recurrence-b-out}%
\begin{tikzpicture}[cheby1scale]
  \fill[contract] (-3,0) -- (-2,0) -- (-2,-1);
  \node[font=\footnotesize] at (3, 0) {$\cyc(n+1)$};
  \node[font=\footnotesize] at (3, -1) {$\cyc(n)$};

  \foreach \v/\x in {k-1/-3, k/-2, k+1/-1, 1/1, 2/2}
    \node[newvert=\v] (t\v) at (\x, 0) {};

  \foreach \v/\x in {k-1/-2, k/-1, 1/1, 2/2}
    \node[newvert=\v] (b\v) at (\x, -1) {};

  \draw[dashed] (tk-1) -- (tk+1);
  \draw[thick] (tk+1) -- (t2);

  \draw[dashed] (bk-1) -- (bk);
  \draw[thick] (bk) -- (b2);
\end{tikzpicture}

\endpgfgraphicnamed
  \caption{Case (b) for the $T\sbspt{n+1}(x)$ recurrence. If the edge
    $(1,2)$ in $\cyc(n)$ is covered by a $t$-path, say the path extends
    back to vertex $k$. The ``long edges'' from $k$ to $1$, and from
    $k+1$ to $1$, represent $n+1-k$ edges in a $t$-path. Expand vertex
    $k-1$ into vertices $k-1$ and $k$, with an edge not in a $t$-path
    between them, and relabel the vertices between the old vertex $k$
    and vertex $1$. Note that $k$ could be $1$, in which case $k-1$ is
    $n$ and no relabeling is necessary. The weight of the new
    configuration is $x$ times the weight of the old.}
  \label{fig:cheby1-t-recurrence-b}
\end{figure}
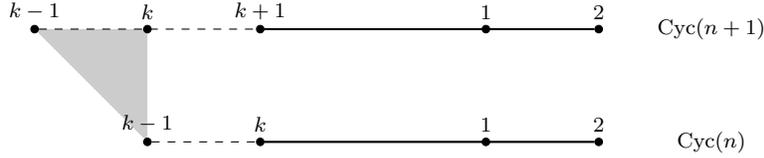

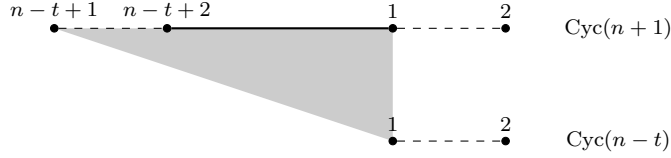
\begin{figure}
  \centering
  \beginpgfgraphicnamed{tikz/cheby1-recurrence-c-out}%
\begin{tikzpicture}[cheby1scale]
  \fill[contract] (-2,0) -- (1,0) -- (1,-1);
  \node[font=\footnotesize] at (3, 0) {$\cyc(n+1)$};
  \node[font=\footnotesize] at (3, -1) {$\cyc(n-t)$};

  \foreach \v/\x in {n-t+1/-2, n-t+2/-1, 1/1, 2/2}
    \node[newvert=\v] (t\v) at (\x, 0) {};

  \foreach \v/\x in {1/1, 2/2}
    \node[newvert=\v] (b\v) at (\x, -1) {};

  \draw[dashed] (tn-t+1) -- (tn-t+2);
  \draw[thick] (tn-t+2) -- (t1);
  \draw[dashed] (t1) -- (t2);

  \draw[dashed] (b1) -- (b2);
\end{tikzpicture}

\endpgfgraphicnamed
  \caption{Case (c) of the $T\sbspt{n+1}(x)$ recurrence. If the $(1,2)$
    edge of $\cyc(n-t)$ is not in a $t$-path, insert a $t$-path preceded
    by an edge immediately behind vertex $1$. The old $(n-t,1)$ edge is
    now an $(n-t, n-t+1)$ edge. This adds a $t$-path and multiplies the
    weight by $-1$.}
  \label{fig:cheby1-t-recurrence-c}
\end{figure}

Now consider $T\sbspt{n-t}(x)$. Case (c) is when the $(1, 2)$ edge is
not in a $t$-path. Insert $t+1$ new edges immediately behind $1$: an
edge not in a $t$-path, followed by a $t$-path. See
\autoref{fig:cheby1-t-recurrence-c}.

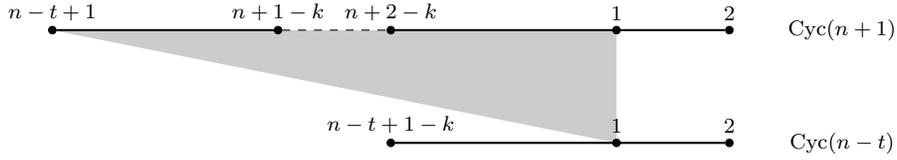
\begin{figure}
  \centering
  \beginpgfgraphicnamed{tikz/cheby1-recurrence-d-out}%
\begin{tikzpicture}[cheby1scale]
  \fill[contract] (-4,0) -- (1,0) -- (1,-1);
  \node[font=\footnotesize] at (3, 0) {$\cyc(n+1)$};
  \node[font=\footnotesize] at (3, -1) {$\cyc(n-t)$};

  \foreach \v/\x in {n-t+1/-4, n+1-k/-2, n+2-k/-1, 1/1, 2/2}
    \node[newvert=\v] (t\v) at (\x, 0) {};

  \foreach \v/\x in {n-t+1-k/-1, 1/1, 2/2}
    \node[newvert=\v] (b\v) at (\x, -1) {};

  \draw[thick] (tn-t+1) -- (tn+1-k);
  \draw[dashed] (tn+1-k) -- (tn+2-k);
  \draw[thick] (tn+2-k) -- (t2);

  \draw[thick] (bn-t+1-k) -- (b2);
\end{tikzpicture}

\endpgfgraphicnamed
  \caption{Case (d) of the $T\sbspt{n+1}(x)$ recurrence. The bottom is
    $\cyc(n-t)$, and edge $(1,2)$ is covered by a $t$-path with $k$
    edges preceding $1$, with $0 \le k \le t-1$. Expand vertex $1$ into
    a sequence of $t+1$ edges of the form ``$t-k$ edges in a $t$-path,
    one edge not in a $t$-path, $k$ edges in a $t$-path''.}
  \label{fig:cheby1-t-recurrence-d}
\end{figure}

Case (d) is shown in \autoref{fig:cheby1-t-recurrence-d}: if the $(1,
2)$ edge of $\cyc(n-t)$ is in a $t$-path, there are $k$ edges in the
path preceding vertex $1$ for some $k$ with $0 \le k < t$. Immediately
behind vertex $1$, add $t+1$ new edges: $k$~edges in a $t$-path, an edge
not in a $t$-path, and $t-k$ edges in a $t$-path. This splits a single
$t$-path into two $t$-paths.

Both operations yield a configuration with one more $t$-path than we
started with, so we've multiplied the weight by $-1$ and have a
configuration for $T\sbspt{n+1}(x)$.

This is a bijection: in $T\sbspt{n+1}(x)$, consider if $(1, 2)$ is in a
$t$-path. If it isn't, is $(n+1, 1)$ in a $t$-path? If no, then the
configuration came from $T\sbspt{n}(x)$ in case~(a); if yes, it came
from case~(c) and $T\sbspt{n-t}(x)$. If $(1, 2)$ is indeed in a
$t$-path: is the path preceded by one, or more than one, edges not in a
$t$-path? If exactly one, the configuration came from case~(d) and
$T\sbspt{n-t}(x)$, and if more than one, it came from $T\sbspt{n}(x)$ in
case~(b). Every possible configuration in $T\sbspt{n+1}(x)$ is accounted
for, so the map is surjective. Since we always simply expand a vertex
into new edges, and because the above argument shows that the cases are
distinguishable, the map is injective.

\subsection{Orthogonality involution}
\label{sec:orth-invol}

Now we address the moments and orthogonality for higher order Chebyshev
polynomials of the first kind. Instead of the Fuss-Catalan numbers, the
moments of $T\sbspnt(x)$ are what we will call \emph{$(t+1)$-reciprocal
  binomial coefficients}:
\begin{equation}
  \label{eq:3}
  B\sbspt{n} = \binom{(t+1)n}{n}.
\end{equation}
Note the similarity to \eqref{eq:fuss-catalan-num}. We say
``$(t+1)$-reciprocal'' because these coefficients give the number of
ways to choose exactly $1/2$, $1/3$, $1/4$, etc., of the elements in a
set.

Let $\Lt$ be the linear functional whose moments are defined by aerated
$(t+1)$-reciprocal binomial coefficients: $\mu\sbspt{(t+1)n}$ equals
$B\sbspt{n}$ and is zero otherwise. We can easily prove an analogue of
\autoref{thm:cheby2-t-integral}:

\begin{theorem}
  \label{thm:cheby1-t-integral}
  Let $n_{1}, n_{2},\dots,n_{k}$ be nonnegative integers. The integral
  \[
  \Lt\left(\prod_{i=1}^{k} T\sbspt{n_{i}} (x)\right)
  \]
  equals the number of ways to mark exactly $1/(t+1)$ of the vertices in
  $\cyc(n_{1}) \sqcup\cdots\sqcup \cyc(n_{k})$ such that no marked
  vertex is followed by $t$ unmarked vertices.
\end{theorem}

\begin{proof}
  The proof follows the now-familiar mantra: the product is the
  generating function for coverings of $\cyc(n_{1}) \sqcup\cdots\sqcup
  \cyc(n_{k})$ by $t$-edge paths with weight $-1$ and with fixed points
  of weight $x$. Applying $\Lt$ can be interpreted as changing all fixed
  points to have weight $1$ and marking exactly $1/(t+1)$ of them. Now
  we apply a sign-reversing involution to the set of those
  configurations: scan through the cycles, and find the first occurrence
  of a $t$-edge path and turn it into a marked vertex followed by $t$
  unmarked vertices, or vice versa.
\end{proof}

\begin{figure}[ht]
  \centering
  \beginpgfgraphicnamed{tikz/cheby1-t-integral-sri-a-out}%
\begin{tikzpicture}[every label/.style={inner sep=1pt},
  polygon/.style={regular polygon, regular polygon sides=#1, minimum
    size=15mm, draw}]

  \node[name=sq, polygon=4] at (0,0) {};
  \node at (sq.corner 2) [vertex, label=135:$1$] {};
  \node at (sq.corner 1) [unvertex, label=45:$2$] {};
  \node at (sq.corner 4) [vertex,  label=315:$3$] {};
  \node at (sq.corner 3) [unvertex, label=225:$4$] {};

  \node[name=pe, polygon=5] at (2.25,0) {};
  \node at (pe.corner 2) [unvertex, label=135:$1$] {};
  \node at (pe.corner 1) [vertex, label=90:$2$] {};
  \node at (pe.corner 5) [unvertex, label=45:$3$] {};
  \node at (pe.corner 4) [unvertex,  label=315:$4$] {};
  \node at (pe.corner 3) [unvertex, label=225:$5$] {};
  \draw[ultra thick] (pe.corner 5) -- (pe.corner 4) --
                     (pe.corner 3) -- (pe.corner 2);

  \node[name=tr, polygon=3] at (4.5, 0) {};
  \node at (tr.corner 1) [vertex, label=90:$2$] {};
  \node at (tr.corner 2) [vertex, label=225:$1$] {};
  \node at (tr.corner 3) [unvertex, label=315:$3$] {};
\end{tikzpicture}
\endpgfgraphicnamed
  \hfill
  \beginpgfgraphicnamed{tikz/cheby1-t-integral-sri-b-out}%
\begin{tikzpicture}[every label/.style={inner sep=1pt},
  polygon/.style={regular polygon, regular polygon sides=#1, minimum
    size=15mm, draw}]
  \node[name=sq, polygon=4] at (0,0) {};
  \node at (sq.corner 2) [unvertex, label=135:$1$] {};
  \node at (sq.corner 1) [unvertex, label=45:$2$] {};
  \node at (sq.corner 4) [unvertex,  label=315:$3$] {};
  \node at (sq.corner 3) [vertex, label=225:$4$] {};

  \node[name=pe, polygon=5] at (2.25,0) {};
  \node at (pe.corner 2) [unvertex, label=135:$1$] {};
  \node at (pe.corner 1) [unvertex, label=90:$2$] {};
  \node at (pe.corner 5) [unvertex, label=45:$3$] {};
  \node at (pe.corner 4) [unvertex,  label=315:$4$] {};
  \node at (pe.corner 3) [unvertex, label=225:$5$] {};
  \draw[ultra thick] (pe.corner 1) -- (pe.corner 2) --
    (pe.corner 3) -- (pe.corner 4);

  \node[name=tr, polygon=3] at (4.5, 0) {};
  \node at (tr.corner 1) [unvertex, label=90:$2$] {};
  \node at (tr.corner 2) [unvertex, label=225:$1$] {};
  \node at (tr.corner 3) [unvertex, label=315:$3$] {};
\end{tikzpicture}
\endpgfgraphicnamed
  \caption{Two (unrelated) configurations in $\L_{3} \big(
    T\sbsp{4}{3}(x) T\sbsp{5}{3}(x) T\sbsp{3}{3}(x)\big)$. Marked
    vertices are indicated by a black circle, unmarked vertices are
    regular corners of the shapes, and the thick edges represent paths
    with $t=3$ edges. In the left configuration, the sign-reversing
    involution would remove the edge in the $5$-cycle and leave the
    $5$-cycle with vertices $2$ and $3$ marked and the others unmarked.
    In the right configuration, vertex $4$ in the $4$-cycle is followed
    by three unmarked vertices, so the involution would replace the
    marked vertex by a path on vertices $4$--$1$--$2$--$3$. Observe that
    it is perfectly acceptable to have no (or all) marked vertices in a
    cycle.}
  \label{fig:cheby1-t-integral-sri}
\end{figure}
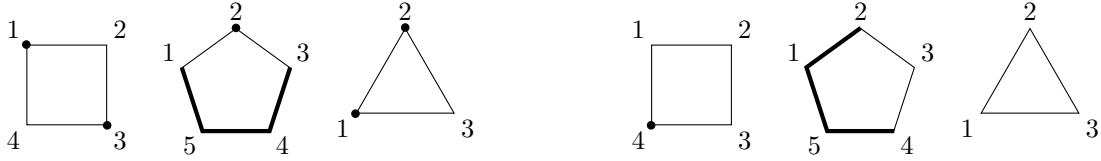

See \autoref{fig:cheby1-t-integral-sri} for two examples of the
involution.

The above theorem immediately gives us the orthogonality relation and
$L^{2}$ norm for $T\sbspnt(x)$.

\begin{corollary}
  \label{thm:cheby1-orthog-l2}
  The Chebyshev polynomials of the first kind and order $t$ are
  $t$-orthogonal with respect to the moments given by aerated
  $(t+1)$-reciprocal binomial coefficients. That is, with $\Lt$ defined
  as above, whenever $m > nt$, then
  \begin{equation}
    \label{eq:cheby1-orthog-l2}
    \Lt \big( U\sbspt{m}(x) U\sbspt{n}(x) \big) = 0
    \quad \text{and} \quad
    \Lt \big( U\sbspt{nt}(x) U\sbspt{n}(x) \big) = t+1.
  \end{equation}
\end{corollary}

\begin{proof}
  Assume $m > nt$. The above theorem tells us that, to find $\Lt \big(
  U\sbspt{m}(x) U\sbspt{n}(x) \big)$, we need to count the number of
  ways to mark exactly $(n+m)/(t+1)$ vertices in $\cyc(m) \sqcup
  \cyc(n)$ such that no marked vertex is followed by $t$ unmarked
  vertices. Since $(n+m)/(t+1)$ is greater than $n$, we must mark at
  least one vertex in $\cyc(m)$, and to insure that we leave no marked
  vertex in $\cyc(m)$ followed by $t$ unmarked vertices, we must mark
  more than $m/t$ vertices there---but $(n+m)/(t+1)$ is strictly smaller
  than $m/t$, so there are zero configurations meeting the criteria.

  When we consider $\Lt \big( U\sbspt{nt}(x) U\sbspt{n}(x) \big)$, we
  must mark exactly $n$ vertices. We can do that and meet the
  marked-unmarked condition by either marking all vertices of $\cyc(n)$,
  or by marking those vertices whose label is congruent to $0, 1,
  2,\dots$ modulo $t$ in $\cyc(nt)$, for a total of $t+1$ different
  configurations.
\end{proof}

\subsection{Generating functions}
\label{sec:cheby1-t-gfs}

In analogy with $UP$ and $UM$, define
\begin{align*}
  TP(t, x, z) &= \sum_{n \ge 0} T\sbspnt(x) z^{n}
  \intertext{and}
  TM(t, z) &= \sum_{n \ge 0} \mu\sbspnt z^{n}.
\end{align*}
These generating functions are not difficult to derive. For the
polynomials, we decompose $T\sbspnt(x)$ by considering vertex $1$. That
vertex is either a fixed point of weight $x$, or is one of the $t+1$
vertices in a $t$-path. If one removes the ``component'' that vertex $1$
is in, the result is a covering of a path. Therefore for $n \ge 1$,
\[
T\sbspnt(x) = x U\sbspt{n-1}(x) - (t+1) U\sbspt{n-(t+1)}(x),
\]
where we take polynomials with negative indices to equal zero. By
multiplying the above equation by $z^{n}$ and summing over $n$, the
above equation immediately yields the generating function for
$T\sbspnt(x)$:
\begin{align}
  \label{eq:cheby1-poly-gf}
  TP(t,x,z) &= 1 + xz UP(t,x,z) - (t+1)z^{t+1} UP(t,x,z) \notag\\
  &= \frac{1 - t z^{t+1}}{1 - (xz - z^{t+1})}.
\end{align}
The generating function for the moments is also similar to that for the
Chebyshevs of the second kind. The recurrence coefficients of
\eqref{eq:cheby1-recurrence} tell us the L-fraction expression
\cite[chapter V, \S 6]{viennot:theorie}:
\begin{equation}
  \label{eq:cheby1-mu-gf-contfrac}
  TM(t, z) = \cfrac{1}{
    1 - \cfrac{(t+1)z^{t+1}}{
      \mylr{1 - \cfrac{z^{t+1}}{
          \mylr{1 - \cfrac{z^{t+1}}{(1 - \cdots)^{t}}}^{t}}}^{t}}}.
\end{equation}
The only recurrence coefficient here that is different from those for
the Chebyshev polynomials of the second kind is the very first one, so
reasoning as in \autoref{fig:dyck-path-decomp}, we have the following:
\begin{equation}
  \label{eq:cheby1-mu-func-eqn}
  TM(t, z) = 1 + (t+1)z^{t+1} (UM(t, z))^{t} TM(t, z).
\end{equation}
That equation also follows from \eqref{eq:cheby1-mu-gf-contfrac} and
the definition of $UM$. The recurrence coefficients also directly tell
us that the moments are the generating function for weighted \luk{}
paths with upsteps and steps $(1,-t)$ in which all steps have weight
$1$ except for a downstep leaving from height $t$---such a step has
weight $t+1$.

\section{Hermite polynomials}
\label{sec:hermite}

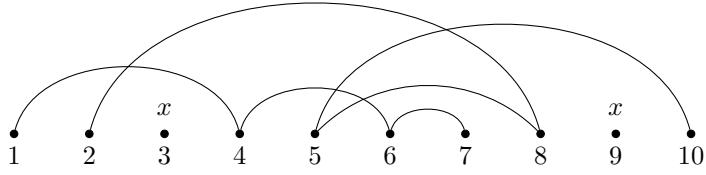
\begin{figure}
  \centering
  \beginpgfgraphicnamed{tikz/hermite-t-poly-example-out}%
\begin{tikzpicture}[edgestyle/.style={looseness=1}]
  \foreach \a in {1,...,10}{
    \node (\a) at (\a,0) [vertex, label=below:$\a$] {} ;
  }

  \drawedge{1}{4} \drawedge{4}{6} \drawedge{6}{7}

  \drawedge{2}{8} \drawedge[out=45, in=135]{5}{8} \drawedge{5}{10}

  \node at (3) [label=above:$x$] {};
  \node at (9) [label=above:$x$] {};
\end{tikzpicture}
\endpgfgraphicnamed
  \caption{A covering of $K_{10}$ by $3$-paths that contributes weight
    $(-1)^{2}x^{2}$ to $H\sbsp{10}{3}(x)$, which equals $x^{10} -
    2520x^{6} + 226800x^{2}$. For clarity, only the edges in covering
    paths are drawn; the remaining $39$ edges of the underlying graph
    are considered understood.}
  \label{fig:hermite-t-poly-example}
\end{figure}

Having addressed paths and cycles, we turn now to complete graphs. The
higher-order matching polynomial for coverings of $K_{n}$ by $t$-paths
is $H\sbspnt(x)$, the Hermite polynomial of order~$t$. These polynomials
are very similar to Chebyshev polynomials of the second kind, but now,
since the underlying graph is the complete graph instead of just a path,
we may have crossings. We will draw coverings for $H\sbspnt(x)$ similar
to how we drew coverings for $U\sbspnt(x)$ but now edges may be between
any two vertices. \autoref{fig:hermite-t-poly-example} shows an example
configuration.

The recurrence relation is hardly any more difficult than that for
higher-order Chebyshev polynomials of the second kind \eqref{eq:cheby2nd-recurrence}. First, observe
that there are $(t+1)!/2$ possible $t$-paths on a set of $t+1$ vertices.
Given $H\sbspt{n+1}(x)$, consider vertex $n+1$: it is either fixed, or
is in a $t$-path with $t$ other vertices, and we have
\begin{equation}
  \label{eq:hermite-t-recurrence}
  H\sbspt{n+1}(x) = x H\sbspnt(x) -
  \binom{n}{t} \frac{(t+1)!}{2} H\sbspt{n-t}(x)
\end{equation}
for $n > 0$, with the usual convention that polynomials with negative
indices are zero and the zeroth polynomial equals one. See Araujo et al.
\cite[\S 2.3]{araujo.estrada.ea:higher-order} for an explicit formula.

For the higher-order Chebyshev $U$ polynomials, the moments were
noncrossing set partitions with all blocks of size $t+1$; now, the
moments involve the same sort of set partitions, but with crossings
allowed. In both cases, the moments are ``complete'' configurations. Let
$\mu\sbspnt$ equal the number of ways to completely cover $K_{n}$ by
$t$-paths and $\Lt$ the corresponding linear functional. Clearly
$\mu\sbspnt$ is zero if $n$ is not a multiple of $t+1$, and if it is a
multiple of $t+1$, then we can count the number of such coverings by
first finding the number of set partitions of $[n]$ with all blocks of
size $t$ (which is just a multinomial coefficient) and then multiplying
by an appropriate power of $(t+1)!/2$, the number of $t$-paths through a
given set of vertices. That is,
\begin{equation}
  \label{eq:hermite-t-mu}
  \mu\sbspt{(t+1)n} = \binom{(t+1)n}{t+1, t+1,\dots,t+1} \left(
    \frac{(t+1)!}{2} \right)^{n}
\end{equation}
with $n$ copies of $t+1$ in the ``denominator'' of the multinomial
coefficient. \nonthmref{thm:hermite-t-orthog} will show that these
numbers really are the moments for the higher-order Hermite polynomials.

For $t=1$, the moments have the very well-known integral representation
as the moments of a positive measure on the real axis, namely
\[
\mu\sbsp{n}{1} = \frac{1}{\sqrt{2 \pi}} \int_{\R} x^{n}
\exp\left(-x^{2}/2\right) \dx.
\]
For $t=2$ and $t=3$, integral representations are also known:
\begin{equation}
  \label{eq:hermite-2-mu-integral}
  \mu\sbsp{3n}{2} = \frac{3^{n}}{\pi} \int_{0}^{\infty}
  x^{n} \sqrt{\frac{2}{3x}} K_{1/3}\left(\frac{2 \sqrt{2x}}{3}\right) \dx,
\end{equation}
where $K_{1/3}$ is the modified Bessel function of the second kind.
This expression follows from the formula given in sequence A025035 of
the OEIS \cite{oeis}. For $t=3$, a formula in sequence A025036 gives a
representation for $\mu\sbsp{4n}{3}$:
\begin{equation}
  \label{eq:hermite-3-mu-integral}
  12^{n}  \int_{0}^{\infty} x^{n-1/4} \frac{3^{1/4}}{2^{1/2}}
  \left(
    \frac{3^{1/2} F(5/4, 3/2) \Gamma(3/4)}{2^{3/4} \pi} -
    \frac{3^{1/4} F(5/4, 3/4)}{\pi^{1/2} x^{1/4}} +
    \frac{F(1/2, 3/4)}{2^{3/4} x^{1/2} \Gamma(3/4)}
  \right) \dx,
\end{equation}
where $F(a, b) = {}_{0}F_{2}\left(\begin{matrix}- \\a & b\end{matrix}\ ;
  -3x/32\right)$, a generalized hypergeometric function.

With $\Lt$ in hand, our next task is to show what happens when one
integrates a product of higher-order Hermite polynomials; like
\autoref{thm:cheby2-t-integral}, this is a generalization of a theorem
of de~Sainte-Catherine and Viennot
\cite[Theorem~2]{sainte-catherine.viennot:combinatorial}.

\begin{theorem}
  \label{thm:hermite-t-integral}
  Let $n_{1}, n_{2},\dots,n_{k}$ be nonnegative integers. The integral
  \[
  \Lt\left(\prod_{i=1}^{k} H\sbspt{n_{i}} (x)\right)
  \]
  equals the number of inhomogeneous coverings of $[n_{1}] \sqcup
  \dots\sqcup [n_{k}]$ by $t$-paths.
\end{theorem}

\begin{proof}
  Just as in \autoref{thm:cheby2-t-integral}, the integral of the
  product is the generating function for complete coverings of $[n_{1}]
  \sqcup \dots\sqcup [n_{k}]$ by two kinds of $t$-paths: black paths,
  which have weight~$-1$ and must stay within on of the $[n_{i}]$, and
  dashed paths, which have weight~$1$ and can go anywhere. We can cancel
  all the black paths with a simple sign-reversing involution: given a
  homogeneous path in $[n_{i}]$, label the path with $(i, j)$, where $j$
  is the smallest index among the $t+1$ vertices of the path. Order
  those labels lexicographically, and switch the first path in that
  ordering from dashed to black, or vice versa. This cancels every
  configuration with a homogeneous $t$-path.
\end{proof}

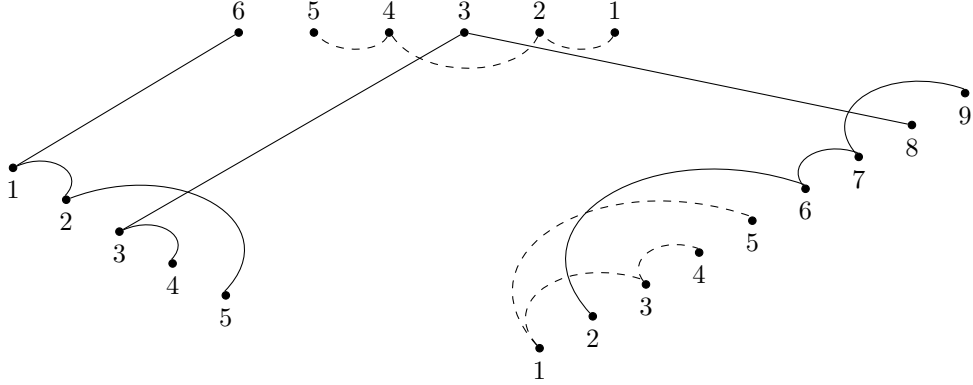
\begin{figure}
  \centering
  \beginpgfgraphicnamed{tikz/hermite-t-integral-out}%
\begin{tikzpicture}[yscale=.6]

\begin{scope}[shift={(-4,4)}]
  \begin{scope}[rotate around={-45:(1,0)}]
    \foreach \a in {1,...,5}{ \node (5-\a) at (\a,0) [vertex,
      label=below:$\a$] {} ; }

    \drawedge{5-1}{5-2} \drawedge{5-2}{5-5}

    \drawedge{5-3}{5-4}
  \end{scope}
\end{scope}

\begin{scope}[shift={(3,0)}]
  \begin{scope}[rotate around={45:(1, 0)}]
    \foreach \a in {1,...,9}{ \node (9-\a) at (\a,0) [vertex,
      label=below:$\a$] {} ; }

    \drawedge[dashed]{9-1}{9-3} \drawedge[dashed]{9-3}{9-4}
    \drawedge[dashed]{9-1}{9-5}

    \drawedge{9-2}{9-6} \drawedge{9-6}{9-7} \drawedge{9-7}{9-9}
  \end{scope}
\end{scope}

\begin{scope}[shift={(4,7)}]
  \begin{scope}[rotate around={180:(1, 0)}]
    \foreach \a in {1,...,6}{ \node (6-\a) at (\a,0) [vertex,
      label=above:$\a$] {} ; }

    \drawedge[dashed]{6-1}{6-2} \drawedge[dashed]{6-2}{6-4}
    \drawedge[dashed]{6-4}{6-5}
  \end{scope}
\end{scope}

\draw (5-1) to (6-6) ;

\draw (5-3) to (6-3) ; \draw (6-3) to (9-8) ;

\end{tikzpicture}
\endpgfgraphicnamed
  \caption{A configuration in $\L\sbsp{}{3}\big(H\sbsp{5}{3}(x)
    H\sbsp{9}{3}(x) H\sbsp{6}{3}(x)\big)$. The sign-reversing involution
    of \autoref{thm:hermite-t-integral} would change the path on
    vertices $1$, $3$, $4$ and $5$ in $[9]$ to black.}
  \label{fig:hermite-t-integral}
\end{figure}

\autoref{fig:hermite-t-integral} shows an example of such a
configuration and the action of the above sign-reversing involution on
it.

The next result, along with uniqueness of moments, implies that the
$\mu\sbspnt$ we defined above are the moments for the higher-order
Hermite polynomials.

\begin{corollary}
  \label{thm:hermite-t-orthog}
  The Hermite polynomials of order $t$ are $t$-orthogonal with respect
  to the moments given by $\mu\sbspnt$: whenever $m > nt$, then
  \begin{equation}
    \label{eq:thm:hermite-t-orthog}
    \Lt \big( H\sbspt{m}(x) H\sbspt{n}(x) \big) = 0
    \quad \text{and} \quad
    \Lt \big( H\sbspt{nt}(x) H\sbspt{n}(x) \big) =
    \left(\frac{t+1}{2}\right)^{n} (nt)!
  \end{equation}
  for $n \ge 0$.
\end{corollary}

\begin{proof}
  The first relation is clear from \autoref{thm:hermite-t-integral},
  since if $m > nt$, there must be a homogeneous $t$-path in $[m]$ and
  so all such configurations are canceled.

  The $L^{2}$ norm requires a bit of work. The involution implies that
  the only configurations we need consider are those in which every
  $t$-path has one vertex in $[n]$, and the other $t$ vertices in
  $[nt]$. To construct such a covering, take a permutation of $[nt]$ and
  insert a bar after every $t$ vertices to form $n$ groups. The first
  group of vertices will be in a path with vertex $1$ from $[n]$, the
  second group in a path with vertex $2$ from $[n]$, and so on. For each
  of the groups, there are $t+1$ ways to insert the vertex from $[n]$
  into the path given by the ordering of the group. However, since the
  paths have no orientation, we have counted every possibility twice, so
  we divide by $2$ for each of the $n$ paths, and obtain the claimed
  $L^{2}$ norm.
\end{proof}

\subsection{Generating functions}
\label{sec:hermite-t-gen-fns}

\begin{figure}
  \centering
  \beginpgfgraphicnamed{tikz/hermite-t-mu-luk-path-decomp-out}%
\begin{tikzpicture}

\upstep{0,0} \wiggle{1,1}{2.5,1}{\delta HM'}
\upstep{2.5,1} \wiggle{3.5,2}{5,2}{\delta^{2} HM'}
\upstep{5,2} \wiggle{6,3}{7.5,3}{\delta^{3} HM'}

\draw (7.5, 3) node[vertex] {} --
  node [midway, above right] {$\lambda\sbsp{t}{t}z$} ++(1, -3)
  node[vertex] {};

\wiggle{8.5,0}{11,0}{HM'}

\end{tikzpicture}
\endpgfgraphicnamed
  \caption{A \luk{} path decomposition for $HM'(t,z)$ for $t=3$. Any
    path contributing to $HM'$ is either empty, or of the above form.
    The $\delta^{k} HM'$ notation means shift all the lower subscripts
    on the $\lambda$'s by $k$.}
  \label{fig:hermite-t-mu-luk-path-decomp}
\end{figure}
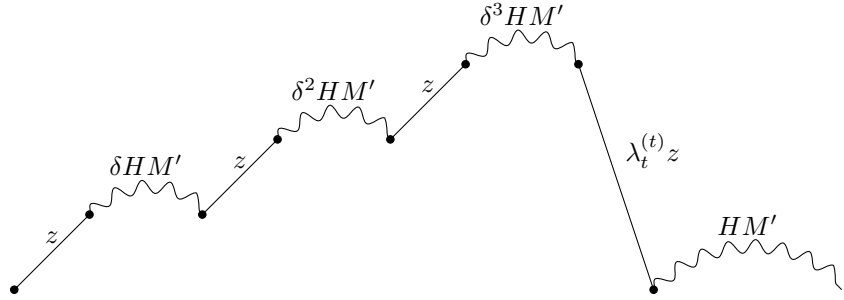

Generating functions for the higher-order Hermite polynomials and their
moments are quite easy to find, since moving from a path to a complete
graph gives ``more symmetry''. First, the polynomials: we define
\[
HP(t, x, z) = \sum_{n \ge 0 } H\sbspnt(x) \frac{z^{n}}{n!}.
\]
Any configuration contributing to $H\sbspnt(x)$ has two kinds of
connected components: $t$-paths and fixed points. The former has weight
$-1$ and there are $(t+1)!/2$ of them on a labeled set of $t+1$ points,
the latter has weight $x$ and there's obviously just one on a point. The
exponential formula \cite[\S 3.3]{aigner:course} immediately gives us
\begin{equation}
  \label{eq:hermite-t-poly-gf}
  HP(t, x, z) = \exp\left(xz - \frac{z^{t+1}}{2}\right).
\end{equation}
The above generating function is a specialization of one found by
Farrell \cite[Theorem~2]{farrell:general}. The same reasoning gives us
the exponential generating function for the moments:
\begin{equation}
  \label{eq:hermite-t-moment-gf}
  HM(t, z) = \sum_{n \ge 0} \mu\sbspnt \frac{z^{n}}{n!} =
  \exp\left(\frac{z^{t+1}}{2}\right).
\end{equation}
We can also find an expression for the ordinary generating function,
since the recurrence coefficients of the polynomials tell us the
continued fraction expansion for that function. We know that the
moments of order $t$ are also the generating functions for \luk{}
paths that consist of upsteps and $t$-downsteps, where upsteps all
have weight $1$ and a $t$-downstep leaving from height $n$ has weight
$\binom{n}{t}(t+1)!/2$ (the recurrence coefficient for
$H\sbspt{n-t}(x)$ in \eqref{eq:hermite-t-recurrence}). For clarity,
let $\lambda\sbspnt = \binom{n}{t}(t+1)!/2$, and write $HM'(t,z)$ for
the ordinary generating function for $\mu\sbspnt$; then, decomposing
\luk{} paths makes the following expression clear:
\begin{equation}
  \label{eq:hermite-t-moment-cf-gf}
  HM'(t, z) =
  \cfrac{1}
    {1 - \cfrac{\lambda\sbspt{t} z^{t+1}}
      {\displaystyle\prod_{k=1}^{t}\mylr{1 - \cfrac{\lambda\sbspt{t+k}z^{t+1}}
        {\displaystyle\prod_{j=1}^{t}\mylr{1 -
            \cfrac{\lambda\sbspt{t+k+j}z^{t+1}}
              {\prod(1 - \cdots)}}}}}}
\end{equation}
This expression, like the continued fraction for $UM(t, z)$ and $TM(t,
z)$, is an L-fraction \cite[\S 6]{viennot:theorie}.
\autoref{fig:hermite-t-mu-luk-path-decomp}, the higher-order Hermite
version of \autoref{fig:dyck-path-decomp}, explains the first steps of
$HM'(3, z)$. Analogous to \eqref{eq:cheby2-mu-functional-eq} and
\eqref{eq:cheby1-mu-func-eqn}, $HM'$ satisfies the functional equation
\[
HM'(t,z) = 1 + \lambda\sbspt{t} z^{t+1} HM'(t,z)
               \prod_{k=1}^{t}\left(\delta^{k}HM'(t,z)\right).
\]
The $\delta$ operator seen here and in
\autoref{fig:hermite-t-mu-luk-path-decomp} is taken from \cite[chap. V,
\S 1]{viennot:theorie}: the generating function $HM'(t, z)$ depends on
$t$, $z$, and the $\lambda\sbspnt$'s, and we could make the dependence
explicit by writing $HM'(t, z, \lambda\sbspt{t},
\lambda\sbspt{t+1},\dots)$. The $\delta$ operator simply increases the
subscript on all the $\lambda$'s: $\delta HM' = HM'(t, z,
\lambda\sbspt{t+1}, \lambda\sbspt{t+2},\dots)$.

\section{Laguerre polynomials}
\label{sec:laguerre}

The final class of higher-order matching polynomials we will consider
are the Laguerre polynomials. The usual Laguerre polynomial
$L_{n}(x^{2})$ is the matching polynomial for the complete bipartite
graph $K_{n,n}$; often the Laguerre polynomials are defined with a
parameter $\alpha$, and while that parameter has combinatorial
meaning---in some sense, it counts cycles; see Foata and Strehl
\cite{foata.strehl:combinatorics}, Labelle and Yeh
\cite{labelle.yeh:combinatorics}, and Simion and Stanton
\cite{simion.stanton:specializations,simion.stanton:octabasic}---we will
not use it; our polynomials correspond to $\alpha=0$.

Let $M_{t}(K_{n,n})$ be the higher-order matching polynomial for
complete bipartite graphs with our usual weights. The degrees of the
polynomials of this sequence are all even, but a $d$-orthogonal sequence
of polynomials requires a polynomial of degree $n$ for every nonnegative
$n$. If $t$ is odd (so that the number of vertices in a $t$-path is
even) we can simply substitute $\sqrt{x}$ for $x$ and get a proper
sequence of polynomials, but if $t$ is even, the resulting matching
polynomials contain both even and odd powers of $x$ and a simple
substitution will not work. We could not find a combinatorially
satisfactory way of converting the matching polynomials for even $t$
into a proper sequence of polynomials, so in this section, we will
hereafter assume that $t$ is an odd positive integer, and define $k$ by
$t = 2k-1$; $k$ is the number of vertices that a $t$-path occupies on
each ``side'' of $K_{n,n}$.

We therefore define the Laguerre polynomial of order $t$ (for odd $t$
only) by the relation
\begin{equation}
  \label{eq:lag-t-poly-defn}
  L\sbspnt(x^{2}) = M_{t}(K_{n,n}).
\end{equation}
Araujo et al. \cite[\S 2.4]{araujo.estrada.ea:higher-order} find
explicit formulas\footnote{Note that the terms of the sum in their
  equation (2.27) are missing a factor of $(-1)^{k}$.} for these
polynomials and derive recurrence relations for the matching
polynomials---but some of the polynomials in those relations are from
graphs not of the form $K_{n,n}$, so we need to derive an appropriate
recurrence relation for $L\sbspnt(x)$ directly. Before we do that,
observe that the number of ways to cover $K_{k,k}$ by a $t$-path is
$(k!)^{2}$: we can orient the path by considering it to start in the
left set of vertices, and then we order the vertices on the left and
right.

\begin{theorem}
  \label{thm:lag-t-recursion}
  With the convention that $L\sbspnt(x) = 0$ when $n < 0$ and
  $L\sbspt{0}(x) = 1$, then for $n \ge 0$, the Laguerre polynomials of
  order $t$ satisfy the recurrence relation
  \begin{multline*}
    L\sbspt{n+1}(x) = x L\sbspt{n}(x) -
     (k!)^{2} \left( \binom{n}{k-1}^{2} + 2 \binom{n}{k} \binom{n}{k-1} \right)
    L\sbspt{n-(k-1)}(x) -{}\\
     \left( \binom{n}{t}\binom{t}{k}
      (k!)^{2} \right)^{2} L\sbspt{n-t}(x).
  \end{multline*}
\end{theorem}

\begin{proof}
  Any covering of $K_{n+1,n+1}$ by $t$-paths can be obtained in one or
  more of the following ways:
  \begin{itemize}
  \item[(a)] Take any covering of $K_{n,n}$ and add two fixed vertices
    at the ``bottom'' of each vertex set; this corresponds to $x
    L\sbspnt$.
  \item[(b)] Take $K_{n+1,n+1}$ and choose $k-1$ vertices among the
    first $n$ vertices of each set. Add in the last vertices of each
    set, put a $t$-path onto those $2k$ vertices, and then ``fill in''
    the rest with any configuration. This corresponds to
    \[
    -(k!)^{2} \binom{n}{k-1}^{2} L\sbspt{n-(k-1)}(x).
    \]
  \item[(c)] Take $K_{n+1,n+1}$ and choose $k$ vertices among the first
    $n$ vertices in the left set, $k-1$ vertices among the first $n$
    vertices in the right set, put a $t$-path down on those $2k$
    vertices, and then fill in the rest with any configuration. This,
    along with exchanging left and right, contributes
    \[
    -2 \binom{n}{k}\binom{n}{k-1}(k!)^{2} L\sbspt{n-(k-1)}(x).
    \]
  \end{itemize}
  Configurations in which vertices $n+1$ on the left and right are both
  in a path, and are not in the same path, are counted twice by item (c)
  above. So we must correct for this by subtracting the total weight of
  those configurations. We need to choose $t$ vertices on each side, and
  then choose $k$ of those vertices to get connected to the bottom
  vertex on the opposite side. Finally, put down two $t$-paths. The
  total contribution of these configurations is
  \[
  \left( \binom{n}{t} \binom{t}{k}(k!)^{2} \right)^{2} L\sbspt{n-t}(x).
  \]

  Adding together the above expressions yields the recurrence relations
  of the theorem.
\end{proof}

\ifhavetikz
  \tikzstyle{mycirc}=[circle, draw=black, inner sep=2]
  \tikzstyle{mycirc}=[circle, draw=black, inner sep=2]
  \tikzstyle{myscaling}=[scale=.6, yscale=-1]
  \newcommand{\mybrace}[4]{\draw[decorate,
    decoration={brace, amplitude=3mm, raise=3mm}] (#1) --
    node[midway, #4=6mm, text width=6em, text centered] {#3} (#2);}
  \newcommand{\mylbrace}[1]{\mybrace{L6}{L1}{#1}{left}}
  \newcommand{\myrbrace}[1]{\mybrace{R1}{R6}{#1}{right}}
\fi

\begin{figure}
  \centering
  \begin{tabular}{cc}
  \subfigure[]{\beginpgfgraphicnamed{tikz/lag-t-recurrence-a-out}%
\begin{tikzpicture}[myscaling]
  \foreach \a in {1,...,6}{
     \node (L\a) at (0, \a) [vertex] {};
     \node (R\a) at (2, \a) [vertex] {};}
  \node (L7) at (0, 7) [mycirc, label=left:fixed] {};
  \node (R7) at (2, 7) [mycirc, label=right:fixed] {};
\end{tikzpicture}
\endpgfgraphicnamed} &
  \subfigure[]{\beginpgfgraphicnamed{tikz/lag-t-recurrence-b-out}%
\begin{tikzpicture}[myscaling]
  \foreach \a in {1,3,4,6} {\node (L\a) at (0, \a) [vertex] {};}
  \foreach \a in {2,5,7} {\node (L\a) at (0, \a) [mycirc] {};}
  \foreach \a in {1,2,3,5} {\node (R\a) at (2, \a) [vertex] {};}
  \foreach \a in {4,6,7} {\node (R\a) at (2, \a) [mycirc] {};}
  \mylbrace{choose $k-1$ vertices}
  \myrbrace{choose $k-1$ vertices}
  \draw[dashed] (L7) -- (R7) -- ++(225:2.25);
  \draw[dashed] (L7) -- ++(-57:2.5);
\end{tikzpicture}
\endpgfgraphicnamed}\\[20pt]
  \subfigure[]{\beginpgfgraphicnamed{tikz/lag-t-recurrence-c-out}%
\begin{tikzpicture}[myscaling]
  \foreach \a in {2,3,4,7} {\node (L\a) at (0, \a) [vertex] {};}
  \foreach \a in {1,5,6} {\node (L\a) at (0, \a) [mycirc] {};}
  \foreach \a in {1,4,5,6} {\node (R\a) at (2, \a) [vertex] {};}
  \foreach \a in {2,3,7} {\node (R\a) at (2, \a) [mycirc] {};}
  \mylbrace{choose $k$ vertices}
  \myrbrace{choose $k-1$ vertices}
  \draw[dashed] (R7) -- ++(235:2.25);
\end{tikzpicture}
\endpgfgraphicnamed} &
  \subfigure[]{\beginpgfgraphicnamed{tikz/lag-t-recurrence-d-out}%
\begin{tikzpicture}[myscaling]
  \foreach \a in {2} { \node (L\a) at (0, \a) [vertex] {}; }
  \foreach \a in {3,5,7} { \node (L\a) at (0, \a) [mycirc] {}; }
  \foreach \a in {1,4,6} { \node (L\a) at (0, \a) [mycirc,
    fill=black!25] {}; }
  \foreach \a in {6} { \node (R\a) at (2, \a) [vertex] {}; }
  \foreach \a in {3,4,7} { \node (R\a) at (2, \a) [mycirc] {}; }
  \foreach \a in {1,2,5} { \node (R\a) at (2, \a) 
    [mycirc, fill=black!25] {}; }
  \mylbrace{choose $t$ vertices \& $k$ among those}
  \myrbrace{choose $t$ vertices \& $k$ among those}
  \draw[dashed] (R7) -- ++(230:2.35);
  \draw[dashed] (L7) -- ++(-63:3);
\end{tikzpicture}
\endpgfgraphicnamed}
  \end{tabular}
  \caption{The four cases of the recurrence relation for
    $L\sbspt{n+1}(x)$; here $t=5$ so $k=3$. In case (a), both bottom
    vertices are fixed; in (b), both bottom vertices are in the same
    path; the bottom vertices are not necessarily adjacent in the
    $t$-path. In (c), the right vertex is in a path, and the bottom
    vertex on the left may or may not be in a path. Case (d) corrects
    the overcounting from (c) when both vertices are in a path: the gray
    circles are in a path together with the bottom vertex from the
    opposite side.}
  \label{fig:lag-t-recurrence}
\end{figure}
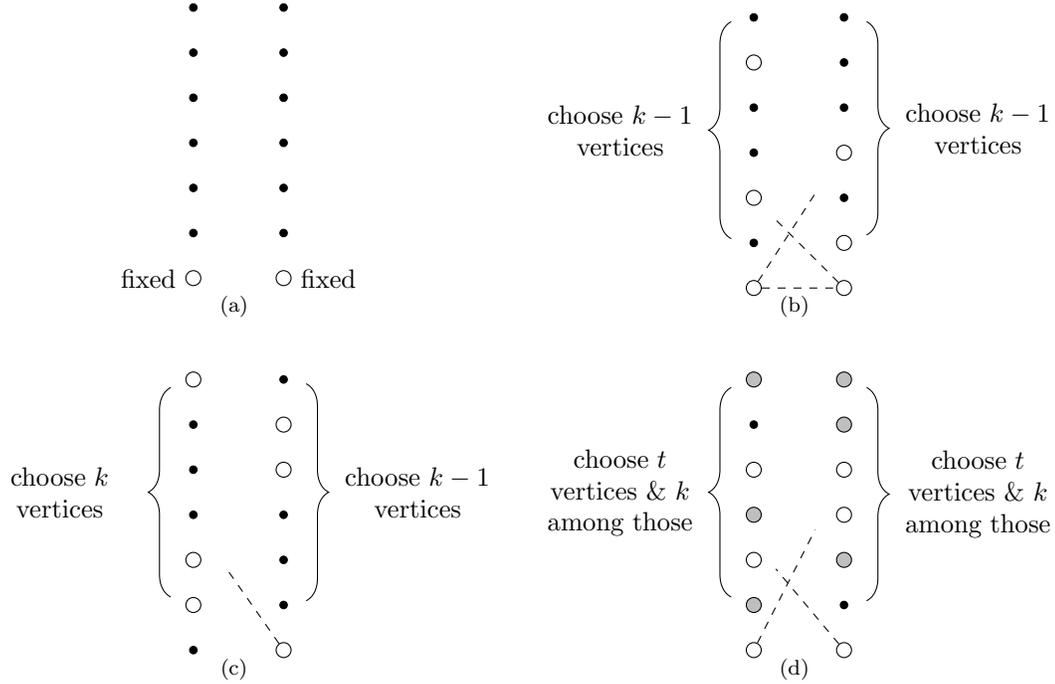

\autoref{fig:lag-t-recurrence} shows the four cases. Observe that if
$k=1$, we indeed recover the recurrence coefficients for the classical
monic Laguerre polynomials: $L_{n+1}(x) = (x - (2n+1))L_{n}(x) - n^{2}
L_{n-1}(x)$.

As usual, we define a linear functional $\Lt$ by $\Lt(x^{n}) =
\mu\sbspnt$, where $\mu\sbspnt$ is the number of complete coverings of
$K_{n,n}$ by $t$-paths, and can count the integral of a product of these
polynomials. This is a generalization of a result of de~Sainte-Catherine
and Viennot \cite[Theorem~5]{sainte-catherine.viennot:combinatorial}.

\begin{theorem}
  \label{thm:laguerre-t-integral}
  Let $n_{1}, n_{2},\dots,n_{j}$ be nonnegative integers. The integral
  \[
  \Lt\left(\prod_{i=1}^{j} L\sbspt{n_{i}} (x)\right)
  \]
  equals the number of inhomogeneous coverings of $K_{n_{1},n_{1}}
  \sqcup \dots\sqcup K_{n_{j},n_{j}}$ by $t$-paths.
\end{theorem}

\begin{proof}
  The proof is the same that we've seen several times now; we start with
  the product of the $L\sbspt{n_{i}}$ and apply $\Lt$, which gives us
  the generating function for complete coverings of the disjoint union
  of the complete bipartite graphs by $t$-paths, in which homogeneous
  $t$-paths may be black (weight $-1$) or dashed (weight $+1$), and
  inhomogeneous paths always are black (weight $+1$). By choosing, say,
  the smallest $i$ such that $K_{n_{i},n_{i}}$ has a homogeneous path,
  finding the homogeneous path with the smallest left vertex inside that
  subgraph, and changing the color from black to dashed or vice versa,
  we have a sign-reversing involution that cancels any configuration
  with a homogeneous $t$-path.
\end{proof}

This immediately implies the $t$-orthogonality of the higher order
Laguerre polynomials (recall that $k = (t+1)/2$):

\begin{corollary}
  \label{thm:laguerre-t-orthog}
  The Laguerre polynomials of order $t$ are $t$-orthogonal with respect
  to the moments given by $\mu\sbspnt$: whenever $m > nt$, then
  \[
  \Lt \left( L\sbspt{m}(x) L\sbspt{n}(x) \right) = 0,
  \]
  and
  \[
  \Lt \left( L\sbspt{nt}(x) L\sbspt{n}(x) \right) = \prod_{i=0}^{n-1}
  \left( \binom{(n-i)t}{t}\binom{t}{k}(k!)^{2} \right)^{2}
  \]
  for $n \ge 0$.
\end{corollary}

\begin{proof}
  The orthogonality relation is clear from
  \autoref{thm:laguerre-t-integral}, since any such configuration in
  $L\sbspt{m}(x) L\sbspt{n}(x)$ must have a homogeneous $t$-path in
  $L\sbspt{m}(x)$. The $L^{2}$ norm can be calculated as follows: after
  applying the sign-reversing involution of the theorem, the only
  remaining configurations are those with $2n$ paths, each with one
  vertex in $K_{n,n}$ and the remaining vertices in $K_{nt,nt}$.
  Consider vertex $1$ on the left and right in $K_{n,n}$: to choose a
  pair of paths that go through those vertices, choose $t$ vertices
  among the $nt$ vertices on the left and right in in $K_{nt,nt}$. There
  are $\binom{nt}{t}^{2}$ ways to do that. Among those $t$ vertices on
  the left, choose $k$ of them to be in the path that goes through
  vertex $1$ on the right side of $K_{n,n}$; the same applies,
  \emph{mutatis mutandis}, on the other side. There are
  $\binom{t}{k}^{2}$ ways to do this. Now take vertex $1$ on the left,
  the $k-1$ vertices not chosen in the second step, and the $k$ vertices
  chosen on the right, and put a $t$-path on those vertices; do the same
  with the remaining vertices---there are $(k!)^{2}$ ways to do that for
  each path. Altogether we've accounted for the $i=n$ factor in the
  above product; now repeat this procedure with vertex $2$ on the left
  and right in $K_{n,n}$ and the remaining $(n-1)t$ vertices in
  $K_{nt,nt}$, and so on; the total number of uncanceled configurations
  is exactly the product above.
\end{proof}

Using just the ``left side of the above argument'' for the $L^{2}$ norm,
we can derive a formula for the moments. We know that the moments for
the Laguerre polynomials of order $t$ are the number of complete
coverings of $K_{n,n}$ by $t$-paths; if $n$ is not a multiple of $k$,
there are zero such coverings, and otherwise if $n = mk$, the number of
coverings is
\begin{equation}
  \label{eq:laguerre-t-moments}
  \mu\sbspt{mk} = \prod_{i=0}^{m-1} \binom{(m-i)k}{k}
  \binom{(m-i)k-1}{k-1} (k!)^{2}.
\end{equation}

\subsection{Generating functions}
\label{sec:lag-t-gf}

The recurrence coefficients found in \autoref{thm:lag-t-recursion} allow
us to give a continued fraction expression for the moment generating
function:
\[
LM(t, z) = \sum_{n \ge 0} \mu\sbspt{n} z^{n}.
\]
The recurrence relation has coefficients in front of
$L\sbspt{n-(k-1)}(x)$ and $L\sbspt{n-t}(x)$, which tells us that the
weighted \luk{} paths whose generating function equals that of the
moments have down steps of $(1, -(k-1))$ and $(1, -t)$; by decomposing
the paths as in \autoref{fig:hermite-t-mu-luk-path-decomp}, we can
express $LM(t, z)$ as an L-fraction:
\begin{equation}
  \label{eq:lag-t-mu-gf}
  LM(t, z) = \frac{1}{1 -
    z^{k-1}\lambda\sbspt{k-1,k-1}\prod_{i=1}^{k-1}\delta^{i} LM(t, z) -
    z^{t}\lambda\sbspt{t,t} \prod_{i=1} \delta^{i} LM(t, z)}
\end{equation}
where $\lambda\sbspt{n, m}$ denotes the weight of a downstep leaving
from height $n$ and falling $m$ steps and $\delta$, as in
\nonthmref[section]{sec:hermite-t-gen-fns}, acts on $LM$ by increasing
the first ``coordinate'' of the coefficients; it changes
$\lambda\sbspt{n, m}$ into $\lambda\sbspt{n+1,m}$.

\section{Further work}
\label{sec:further-work}

In this paper, we've worked with sets of polynomials that satisfy a
recurrence of order $t$ and in each case, found a linear functional with
respect to which the polynomials are $t$-orthogonal. However, as Van
Iseghem \cite{iseghem:approximants} and Maroni
\cite{maroni:lorthogonalite} have shown, such sequences of polynomials
are naturally associated to not just a single linear functional, but
$t-1$ of them. These functionals are defined by
\[
\L\sbspt{k}(P_{n}(x) P_{m}(x)) =
\begin{cases}
  0 & \text{if $n > tm + k$,}\\
  \text{nonzero} & \text{if $n = tm + k$,}
\end{cases}
\]
for $k = 0,\dots,t-2$. The $n$th moment of $\L\sbspt{k}$ is the
generating function for weighted \luk{} paths of length $n$ that end at
height $k$. In this work we have only addressed the $k=0$ functionals
and in light of the results of Van Iseghem and Maroni, the combinatorial
theory of these higher-order matching polynomials is not entirely known
until interpretations of the higher functionals are known.

\renewcommand{\MR}[1]{%
  \href{http://www.ams.org/mathscinet-getitem?mr=#1}{MR~#1}}
\newcommand{\ISBN}[1]{\href{http://worldcat.org/isbn/#1}{ISBN~#1}}
\bibliographystyle{amsplainurl}

\bibliography{hompdo}

\end{document}